\documentclass[11pt,a4paper,leqno]{amsart}
\usepackage{amssymb}
\newif\ifdiags
\ifdiags\usepackage[PostScript=dvips,nohug]{diagrams}\fi
\theoremstyle{plain}
\newtheorem{thm}{Theorem}[section]
\newtheorem{prop}[thm]{Proposition}
\theoremstyle{definition}
\newtheorem{defn}[thm]{Definition}
\newcommand{\emphdef}{\textit}
\numberwithin{equation}{section}

\newcommand{\thismonth}{\ifcase\month\or
  January\or February\or March\or April\or May\or June\or
  July\or August\or September\or October\or November\or December\fi
  \space\number\year}
\newcommand{\bauth}[1]{\mbox{#1},}
\newcommand{\bart}[1]{\textit{#1},}
\newcommand{\bjourn}[3]{#1 \textbf{#2} (#3)}
\newcommand{\bbook}[1]{\textsl{#1},}
\newcommand{\bpp}[2]{pp.~#1--#2.}
\DeclareMathAlphabet{\mathrmsl}{OT1}{cmr}{m}{sl}
\newcommand{\symb}[2]{\newcommand{#1}{\mathit{#2}} }
\newcommand{\rssymb}[2]{\newcommand{#1}{\mathrmsl{#2}} }
\newcommand{\rsoper}[3][n]{\newcommand{#2}{\mathop{\mathrmsl{#3}}%
\ifx n#1\nolimits\else\limits\fi} }
\symb\vol{vol}
\rssymb\iden{id}
\rssymb{\scal}{scal}
\rsoper\divg{div}
\rsoper\trace{tr}
\rsoper\sym{sym}
\newcommand{\R}{{\mathbb R}}
\newcommand{\C}{{\mathbb C}}

\renewcommand{\geq}{\geqslant}
\renewcommand{\leq}{\leqslant}

\newcommand{\tens}{\mathbin{\otimes}}
\newcommand{\dual}{^{*\!}}
\newcommand{\Cinf}{\mathrm{C}^\infty}

\newcommand{\MS}{{\mathcal H}}

\newcommand{\cL}{{\mathcal L}}

\begin{document}
\title{Two dimensional Einstein-Weyl structures}
\author{David M. J. Calderbank}
\address{Department of Mathematics and Statistics\\
University of Edinburgh\\ King's Buildings, Mayfield Road\\
Edinburgh EH9 3JZ\\ Scotland.}
\email{davidmjc@maths.ed.ac.uk}
\date{\thismonth}
\begin{abstract} All local solutions of the two dimensional
Einstein-Weyl equations are found, and related to the compact
examples obtained in~\cite{DMJC1}.
\end{abstract}
\maketitle
\section{Introduction}

Einstein-Weyl geometry has received much attention in recent
years~\cite{CP,Gauduchon4}, particularly in three
dimensions~\cite{Hitchin3,PT}, where Einstein-Weyl structures arise as
symmetry reductions of the self-duality equations for four dimensional
conformal structures~\cite{JT}. An Einstein-Weyl structure on an $n$-manifold
$M$, with $n\geq3$, consists of a conformal structure together with a
compatible (i.e., conformal) torsion-free connection $D$ such that the
symmetric trace-free part of the Ricci tensor of $D$ vanishes. When $D$ is the
Levi-Civita connection of a compatible Riemannian metric then this metric is
Einstein. As with Einstein metrics, the two dimensional story is somewhat
exceptional. A conformal surface with compatible torsion-free connection $D$
is said to be Einstein-Weyl~\cite{DMJC1} iff
$$D\scal^D-2\divg^D F^D=0,$$ where $\divg^D=\trace D$ is the divergence on
$2$-forms, $F^D$ is the Faraday 2-form of $D$, which is the curvature of $D$
on a natural real line bundle $L^1$, and $\scal^D$ is the scalar curvature of
$D$ viewed as a section of $L^{-2}:=(L^1)\dual\tens (L^1)\dual$. If $F^D=0$,
then $D$ is locally the Levi-Civita connection of a metric of constant scalar
curvature.

The idea of studying the two dimensional case was first suggested
in~\cite{PT}, in which Pedersen and Tod proposed the goal of classifying the
compact examples. This classification was carried out
in~\cite{DMJC1}. Pedersen and Tod also claimed that the local solutions should
depend on a single holomorphic function of one variable. The main aim of this
paper is to show that this is true for the definition above and to obtain all
the solutions explicitly in terms of this holomorphic function.
\begin{thm}\label{main} Let $D$ be an Einstein-Weyl structure in two
dimensions. Then there is a local complex coordinate $\zeta=x+iy$ and a
holomorphic function $h$ such that $D=D^g+\omega$, where $g=dx^2+dy^2$ is the
flat metric and
$$\omega=\frac1{\overline h-\zeta}d\zeta
+\frac1{h-\overline\zeta}d\overline\zeta.$$
\end{thm}

The notation used in this paper follows~\cite{DMJC1}. In particular $L^w$ is
the real line bundle associated to the representation $A\mapsto|\det A|^{w/2}$
of $GL_2(\R)$, so that $L^{-2}$ may be identified with $\Lambda^2T\dual M$
once an orientation is chosen. A conformal structure on $M$ may be viewed as a
metric on $TM$ with values in $L^2$. A \emphdef{Weyl derivative} is a
covariant derivative $D$ on $L^1$. Each choice of compatible metric $g$
trivialises $L^1$. If the corresponding trivial Weyl derivative is denoted
$D^g$, then $D=D^g+\omega$ for some connection $1$-form $\omega$. It is
well known that Weyl derivatives on a conformal manifold correspond
bijectively to compatible torsion-free connections. For instance, $D^g$
corresponds to the Levi-Civita connection of $g$.

I prove theorem~\ref{main} in section~\ref{class}. In the following sections I
discuss the extent to which the solutions are genuinely distinct Einstein-Weyl
structures, explain the geometry behind the solutions and show how the compact
examples arise when $h$ is a (possibly degenerate) M\"obius transformation.
I end the paper, with a brief discussion of the ``twistor theory'' of
M\"obius structures.

\section{Local solution of the two dimensional Einstein-Weyl equations}
\label{class}

The two dimensional Einstein-Weyl condition is, a priori, nonlinear, but may
in fact be linearised. In order to do this I shall make use of the
relationship between Weyl structures and M\"obius structures~\cite{DMJC1}.
\begin{defn} A \emphdef{M\"obius structure} on a conformal manifold
$M$ is a (smooth) second order linear differential operator $\MS$ from
$L^1$ to $S^2_0T\dual M\tens L^1$ such that for some Weyl derivative $D$,
the operator $\MS-\sym_0 D^2$ is zero order.
\end{defn}
A M\"obius structure is a possibly non-integrable and unoriented version of a
complex projective structure. More precisely, a M\"obius structure $\MS$
possesses a tensorial invariant $C^\MS\in\Cinf(M,L^{-2}\tens T\dual M)$ called
the \emphdef{Cotton-York tensor} of $\MS$, by analogy with the three
dimensional case. The M\"obius structure is integrable (i.e., given locally by
the trace-free Hessian in a suitable chart) iff $C^\MS=0$ (see~\cite{DMJC1}).
In this case, if $M$ is oriented and $\phi$ is a local orientation preserving
conformal diffeomorphism then $\phi^*\MS-\MS$ can be identified with the
Schwarzian derivative of $\phi$, and so the M\"obius structure defines a
complex projective structure.

In general the Cotton-York tensor of $\MS$ may be computed using an arbitrary
Weyl derivative $D$. The result is:
$$C^\MS=\divg^D\left(r^D_0-\tfrac14\scal^D\iden+\tfrac12F^D\right),$$
where $r^D_0=\MS-\sym_0 D^2$. From this, the following result is immediate.
\begin{prop} A Weyl structure $D$ in two dimensions is Einstein-Weyl
if and only if the trace-free Hessian $\sym_0 D^2$ is locally the trace-free
Hessian in some conformal chart.
\end{prop}

Consequently, if $D$ is Einstein-Weyl, there is locally a flat metric $g$ such
that $\sym_0D^2=\sym_0(D^g)^2$. If $D=D^g+\omega$, then $\sym_0
D^g\omega-\omega\tens_0\omega=0$. Solving this will give all local solutions
of the Einstein-Weyl equation.

Although this equation is still nonlinear, its resemblance to the Riccati
equation suggests a way of linearising it. To do this, let $\zeta$ be a local
complex coordinate such that $g=d\zeta\,d\overline\zeta$ and write $\omega=f
d\zeta+\overline f d\overline\zeta$ for some complex-valued function $f$. Then
the equation for $\omega$ becomes $f'=f^2$, where $f'$ denotes the complex
linear part of $df$. This is the Riccati equation if $f$ is
holomorphic. Substituting $f=-u^{-1}u'$ (which is always possible locally)
gives $u''=0$ and so $u'=\overline h_0$ for some holomorphic function
$h_0$. Hence $u=\overline h_0(\zeta-\overline h)$, where $h$ is also
holomorphic, and so $f=-u^{-1}u'=1/(\overline h-\zeta)$.

This proves Theorem~\ref{main}.

\section{Gauge transformations}

In order to show that the Einstein-Weyl solutions of Theorem~\ref{main} depend
in an essential way on a single holomorphic function, it is necessary to ask
to what extent the solutions are equivalent under a change of complex
coordinate $\zeta$.

An initial observation is that the scalar curvature and Faraday curvature of
$D$ are given by the real and imaginary parts of $h'/(h-\overline\zeta)^2$. In
particular, $D$ is flat if and only if $h$ is constant, and so most of the
solutions are non-trivial.

More generally, note that the complex coordinate $\zeta$ has been partially
fixed by requiring that the trace-free Hessian induced by this coordinate
chart is the M\"obius structure determined by $D$. Hence, the only remaining
freedom in $\zeta$ is the freedom to apply M\"obius transformations.

If $\zeta=\phi(z)=(az+b)/(cz+d)$ with $ad-bc\neq0$, then
$$d\zeta\,d\overline\zeta=\left|\frac{ad-bc}{(cz+d)^2}\right|^2
dz\,d\overline z.$$
After rescaling the metric, the Einstein-Weyl structure is given by the new
holomorphic function $\tilde h=\overline\phi^{-1}\circ h\circ\phi$, where
$\overline\phi(z)=(\overline az+\overline b)/(\overline cz+\overline d)$.
Thus the Einstein-Weyl structure determines $\overline h$ up to conjugation by
a M\"obius transformation.

\section{Geometry of Weyl connections}

The transformation law for $\overline h$ may be traced back to the fact that
it defines a Weyl derivative $D$. If $J^1L^1$ denotes the bundle
of 1-jets of $L^1$, then $D$ is a section of the affine subbundle $A(M)$ of
$L^{-1}\tens J^1L^1$ given by the splittings of the 1-jet projection
$J^1L^1\to L^1$. This affine bundle is modelled on $T\dual M$.

A M\"obius structure on $M$, as a second order linear differential operator
on $L^1$, defines a vector subbundle $E(M)$ of the $2$-jet bundle $J^2L^1$.
Since this operator is given in coordinates by the trace-free Hessian plus
a zero order term, the $1$-jet projection $E(M)\to J^1L^1$ is surjective,
and its kernel, which is the intersection of $E(M)$ with
$S^2T\dual M\tens L^1$, is the line bundle $L^{-1}$, embedded as the
trace-like tensors.

If $\mu$ is a nonvanishing section of $L^1$ and $g$ is the compatible metric
corresponding to this trivialisation, then $\scal^g\mu^2$ is a function whose
value at $x$ depends quadratically on $(j^2\mu)_x$.  This turns out to define
a natural metric of signature $(3,1)$ on $E(M)$ such that the distinguished
line $L^{-1}$ is null and is the only null line in the kernel of the
projection from $E(M)$ to $L^1$ (see~\cite{DMJC1}, and also~\cite{GNotes} for
more details in the analogous higher dimensional case). Consequently, there is
a natural sphere bundle $S^2(M)$ over $M$, namely the space of null lines in
$E(M)$, and this sphere bundle has a distinguished section. The complement of
this section is an affine bundle and this affine bundle is canonically
isomorphic to $A(M)$ by projecting each null line into $J^1L^1$. Therefore a
Weyl connection is a section of $S^2(M)$ which does not meet the distinguished
section.

Now suppose that the M\"obius structure is integrable. Then $E(M)$ also
possesses a canonical flat connection compatible with the Lorentzian
structure. This flat connection identifies $S^2(M)$ locally with $M\times
S^2$, and the distinguished section gives the developing map from (open
subsets of) $M$ to $S^2$. A complex coordinate $\zeta$ on $M$ compatible with
the M\"obius structure identifies this sphere of parallel sections with
$\C\cup\{\infty\}$, so that $\zeta$ itself corresponds to the distinguished
section of $S^2(M)$. The function $\overline h$ arising in Theorem~\ref{main}
is therefore the local coordinate representation of an antiholomorphic section
of $S^2(M)$. The expression $1/(\overline h-\zeta)$ may be viewed as
stereographic projection from $S^2(M)$ onto $A(M)$. It is well defined for
$\overline{h(\zeta)}\neq\zeta$ and sends poles of $h$ to the origin of $A(M)$
determined by the Levi-Civita connection of $g$.

In fact, if the Weyl connection $D$ is viewed as a section of $A(M)$, its
covariant derivative (as a section of $T\dual M\tens V(A(M))=T\dual M\tens
T\dual M$) can be identified with $r^D_0+\frac14\scal^D\iden-\frac12F^D$,
where $r^D_0=\MS-\sym_0D^2$ (cf.~\cite{GNotes}). Hence $D$ is holomorphic iff
it is flat, and antiholomorphic (with respect to $\MS$) iff $\MS=\sym_0D^2$.
The apparent nonlinearity of the Einstein-Weyl condition arises from the fact
that the flat connection on $A(M)$ is not affine. Nevertheless, it identifies
$A(M)$ locally with an open subset of $M\times S^2$, and so the condition for
a section to be antiholomorphic is in fact linear.

\section{The compact examples}

In~\cite{DMJC1}, the local forms of the Einstein-Weyl structures on compact
surfaces were found. In this section I will show that these solutions are
obtained when $h$ is a (possibly degenerate) M\"obius transformation.

The solutions are given explicitly in terms of a compatible metric and
connection $1$-form as follows:
\begin{align*}
g&=P(v)^{-1}dv^2 + v^2 dt^2\\
\omega&=A v^2\,dt,\\
\tag*{\text{where}} P(v) &=-A^2v^4 + Bv^2 +C,
\end{align*}
and $A,B,C$ are arbitrary constants, constrained only by the condition that
$P(v)$ should be somewhere positive. In~\cite{DMJC1}, I showed that these
Einstein-Weyl structures are defined on $S^2$ (for $C>0$) or $T^2$ (for $C<0$)
by writing $v$ as a elliptic function of $x$ so that $v'(x)^2=P(v)$.
If instead, one substitutes $v^2=1/u$ and rescales $g$ and $t$ by $2$, then
the Einstein-Weyl structure becomes
\begin{align*}
g&=\frac1u\left(\frac{du^2}{-A^2+Bu+Cu^2} + dt^2\right)\\
\omega&=\frac{A\,dt}{2u}.
\end{align*}
Now for $C>0$ introduce a new coordinate $r$ by
$u'(r)^2=(-A^2+Bu+Cu^2)/(Cr^2)$. This is readily integrated to give
$$u(r)=\frac{(B^2+4A^2C)-2Br^2+r^4}{4Cr^2}.$$ Rescaling so that the
metric is $dr^2+r^2dt^2$ leads to the solution of Theorem~\ref{main}
given by $h(\zeta)=(B-2iA\sqrt{C})/\zeta$. Notice that
$\overline{h(\zeta)}=\zeta$ iff $\zeta\overline\zeta=(B-2iA\sqrt{C})$. Hence
if $A\sqrt{C}\neq0$, the solution is globally defined on $S^2$.

For $C<0$ introduce instead a coordinate $\theta$ by
$u'(\theta)^2=(-A^2+Bu+Cu^2)/(-C)$. This integrates to give
$$u(\theta)=\frac{B+\sqrt{B^2+4A^2C}\sin\theta}{-2C}$$
and the Einstein-Weyl structure becomes
\begin{align*}
g&=\frac1{B+\sqrt{B^2+4A^2C}\sin\theta}\left(dt^2+d\theta^2\right)\\
\omega&=\frac{A\sqrt{-C}\,dt}{B+\sqrt{B^2+4A^2C}\sin\theta},
\end{align*}
which is globally defined on $T^2$ (for $C<0$ and $B^2+4A^2C>0$). After
rescaling so that the metric is $e^{2t}(dt^2+d\theta^2)$, the solution
$h(\zeta)=i(B+2A\sqrt{-C})\zeta/\sqrt{B^2+4A^2C}$ of Theorem~\ref{main} is
obtained.

More generally if $\overline h$ is an orientation reversing M\"obius
transformation, then the Weyl connection is well defined away from the fixed
points of this transformation. Hence the elliptic elements, apart from the
simple inversions (which have an invariant circle), give solutions globally
defined on $S^2$ (equivalent to one of the solutions above).  The hyperbolic
elements, with two fixed points, correspond to the solutions on $T^2$ (they
are periodic solutions on a cylinder).  The remaining cases occur as
limits. For instance the simple inversions, such as $\zeta\mapsto
1/\overline\zeta$, give the hyperbolic metric.

\section{Twistor theory}

The \emphdef{twistor space} of a conformal $2$-manifold $M$ is its orientation
double cover, viewed as a complex curve $\Sigma$. This is a rather trivial two
dimensional analogue of the four dimensional theory (see, for
instance,~\cite{Gauduchon4}). Note that $\Sigma$ has a real structure given by
the nontrivial involution in each fibre and that $M$ may be recovered from
$\Sigma$ as the moduli space of real pairs of points. The full moduli space of
(unordered, distinct) pairs of points in $\Sigma$ is
$M^\C=\bigl(\Sigma\times\Sigma\smallsetminus\Delta(\Sigma)\bigr)/S_2$.  This
complex surface has a natural conformal structure: a tangent vector to $M^\C$
at $\{x_1,x_2\}$ consists of a pair of tangent vectors to $\Sigma$ (at $x_1$
and $x_2$), and it is null if one of these components vanishes. Hence $\Sigma$
is (locally) the space of null geodesics in $M^\C$. Of course $M^\C$ is the
natural space in which real analytic functions on $M$ may be written
$f=f(z,\overline z)$ with $f$ holomorphic in two variables.

Although this notion of twistor space has no real content, it does provide a
formal way to distinguish an integrable M\"obius structure in two dimensions
from a one dimensional complex projective structure. The former is a
trace-free Hessian $L^1\to S^2_0T\dual M\tens L^1$, whereas the latter is a
second order operator $\cL\to(T^*\Sigma)^2\tens\cL$ (on a line bundle $\cL$
with $\cL^2=T\Sigma$) whose symbol is the identity.  The two are easily
related: $(T^*\Sigma)^2$ is the pullback of $S^2_0T\dual M$, and since
$T\Sigma\tens\overline{T\Sigma}$ is (the pullback of) $L^2\tens\C$, it follows
that $\cL\tens\overline\cL$ can be identified with $L^1\tens\C$. The
projectivisation of $J^1\cL$ corresponds to $S^2(M)$, and the complex
projective structure defines a connection $J^1\cL\to J^2\cL\leq J^1(J^1\cL)$
which projectivises to the flat connection on $S^2(M)$ induced by the
integrable M\"obius structure.

A more satisfying twistorial description would encode the M\"obius structure
in pure holomorphic geometry. Nevertheless, I hope the na\"\i ve twistor
theory given here at least provides some light entertainment.

\end{document}